\documentclass[a4paper,11pt]{article}

\usepackage[english]{babel}
\usepackage{amsmath}
\usepackage{amsfonts}
\usepackage{amssymb}
\usepackage{amsthm}
\usepackage{epsfig,graphics}
\usepackage{graphicx}
\usepackage{hyperref}

\title{Wonderful models for toric arrangements}

\author{Luca Moci\footnote{Institut fuer Mathematik, Technische Universitaet Berlin,
Strasse des 17. Juni, 136, Berlin. moci@math.tu-berlin.de
}}

\begin{document}

\newtheorem{te}{Theorem}[section]
\newtheorem{lem}[te]{Lemma}
\newtheorem{pr}[te]{Proposition}
\newtheorem{prob}[te]{Problem}
\newtheorem{co}[te]{Corollary}
\newtheorem{cj}[te]{Conjecture}
\theoremstyle{definition}
\newtheorem{re}[te]{Remark}
\newtheorem{ex}[te]{Example}

\newcommand{\lc}{\Lambda_C}
\newcommand{\lcb}{\overline{\Lambda_C}}
\newcommand{\cs}{\mathcal{S}}
\newcommand{\s}{\mathcal{S}}
\newcommand{\g}{\mathcal{G}}
\newcommand{\ir}{\mathcal{I}}
\newcommand{\vs}{\widetilde{\mathcal{V_S}}}
\newcommand{\vt}{\widetilde{\mathcal{V_Q}}}
\newcommand{\vsp}{\mathcal{V_S}}
\newcommand{\vtq}{\mathcal{V_Q}}
\newcommand{\rxt}{\mathcal{R}_{\widetilde{X}}}
\newcommand{\pcn}{\mathbb{P}(\mathbf{N}_T(C))}
\newcommand{\pc}{\mathbb{P}_C}
\newcommand{\zxt}{\mathbf{Z}_{\widetilde{X}}}
\newcommand{\zxg}{\mathbf{Z}_{{\widetilde{X}},\g}}
\newcommand{\yxt}{\mathbf{Y}_{\widetilde{X}}}
\newcommand{\zs}{z^{\mathcal{S}}}
\newcommand{\zt}{z^{\mathcal{Q}}}
\newcommand{\dn}{\mathbf{D}_\mathcal{N}}

\maketitle


\begin{abstract}
We build a wonderful model for toric arrangements.
We develop the "toric analogue" of the combinatorics of nested sets, which allows to define a family of smooth open sets covering
the model. In this way we prove that the model is smooth, and we give
a precise geometric and combinatorial description of the normal crossing divisor.
\end{abstract}

\section{Introduction}

In the spirit of the much studied, ingenious construction of wonderful
models for arrangements of linear subspaces by De Concini and Procesi, we here present a construction of wonderful
models for toric arrangements. The latter have attracted much interest in
recent years, since they proved to be deeply related with a wide number of
topics, including vector partition functions and integral points in
polytopes (\cite{DPt}, \cite{li}), matroids and zonotopes (\cite{MoT}), Lie theory (\cite{L2}, \cite{Mcm}, \cite{Mo}),
and index theory (\cite{DPV2}).

Toric arrangements can be viewed as periodic counterparts of linear arrangements.
Accordingly, their theory is much inspired by the theory of
hyperplane arrangements as it evolved over the last decades (see in particular \cite{Dou}, \cite{DPt},
 \cite{MS}). The quest for a wonderful model of toric arrangements is then a natural next step.


Let $T$ be a complex torus and $\Lambda$ its group of characters.
Let $\widetilde{X}$ be a finite subset of $\Lambda\times \mathbb{C}^*$.
For every pair $(\lambda,a)\in \widetilde{X}$ we define the hypersurface of $T$:
$$H_{\lambda,a}\doteq \left\{ t\in T |\lambda(t)-a=0 \right\}.$$
The collection
$$\mathcal{T}_{\widetilde{X}}\doteq\left\{H_{\lambda,a}, (\lambda,a)\in {\widetilde{X}}\right\}$$ is called the \emph{toric arrangement} defined by ${\widetilde{X}}$ on $T$.

Let $\rxt$ be the complement of the arrangement:
$$\rxt\doteq T\setminus\bigcup_{(\lambda,a)\in {\widetilde{X}}}H_{\lambda,a}.$$

In this paper we build a smooth minimal model $\zxt$ containing $\rxt$ as an open set with complement a normal crossing divisor $\mathbf{D}$, and a proper map  $\pi: \zxt \rightarrow T$ extending the identity of $\rxt$. We call $\zxt$ the \emph{wonderful model} of $\mathcal{T}_{\widetilde{X}}$, in analogy to the wonderful model built by De Concini and Procesi \cite{DPw} for arrangements of subspaces in a vector (or projective) space.

The model $\zxt$ has several applications. In particular, it potentially is a powerful tool to describe the cohomology ring $H^*(\rxt, \mathbb{Q})$.
This computation is one of the central, most outstanding questions in the theory of toric arrangements. In the special case of arrangements defined by \emph{unimodular} lists of vectors, it has been answered in \cite{Lo}, \cite{DPt}. We plan to face this problem in a future paper, by applying to the model $\zxt$ the general method described in \cite{M}.
\bigskip

The paper is organized as follows. In Section 2 we give the first definitions, we make some basic remarks and we build the wonderful model. In Section 3 we develop the necessary combinatorial tools, i.e the "toric analogues" of the notions of irreducible set, building set, nested set, and adapted basis. In Section 4 we define some smooth open sets of the model and we prove that they cover $\zxt$. In Section 5 the open sets are used to prove that the complement of $\rxt$ in $\zxt$ is a normal crossing divisor, and to describe its irreducible components and their intersections (see Theorem \ref{ncd}).

\paragraph{Acknowledgements}
I wish to thank my supervisor Corrado De Concini for many illuminating suggestions and helpful remarks. I am also grateful to Maria Angelica Cueto, Jacopo Gandini and Bernd Sturmfels for stimulating discussions.

\section{First definitions and remarks}

\subsection{Toric arrangements}

Let $\Lambda$ be a lattice of rank $n$ and $U=\Lambda\otimes_{\mathbb{Z}} \mathbb{C}$  the complex vector space obtained by extending the scalars of $\Lambda$.

Let ${\widetilde{X}}$ be a finite set in $\Lambda\times \mathbb{C}^*$, and set
$$X\doteq \{\lambda |(\lambda,a)\in {\widetilde{X}}\}.$$

Given $A\subseteq X$, we denote by $\langle A\rangle_{\mathbb{Z}}$  and $\langle A\rangle_{\mathbb{C}}$ respectively the sublattice of $\Lambda$ and the subspace of $U$ spanned by $A$.
We will always assume the sublattice $\langle X\rangle_{\mathbb{Z}}$ to have finite index in $\Lambda$;
otherwise we can replace $\Lambda$ with $\Lambda\cap \langle X\rangle_{\mathbb{C}}$.

Then we define
$$T\doteq {Hom(\Lambda,\mathbb{C}^*)}.$$
The group $T$ is isomorphic to $(\mathbb{C}^*)^n$, and its group of characters $Hom(T, \mathbb{C}^*)$ is identified with $\Lambda$.
Indeed given $\lambda\in \Lambda$
and $t\in T$, we can take any representative $\varphi_t\in
Hom(\Lambda,\mathbb{C})$ of $t$ and set
$$\lambda(t)\doteq e^{2\pi i \varphi_t(\lambda)}.$$

For every pair $(\lambda,a)\in {\widetilde{X}}$ we define:
$$H_{\lambda,a}\doteq \left\{ t\in T |\lambda(t)-a=0 \right\}.$$

We remark that in general the hypersurfaces $H_{\lambda,a}$ are not connected; and even if they are, their intersections are not (see Remark \ref{pri} and Example \ref{eC2} below). Then we consider the set $\mathcal{C}({\widetilde{X}})$ of all the connected components of all the intersections of the hypersurfaces $H_{\lambda,a}$. This is a poset (with respect to inclusion) which plays a major role in the study of toric arrangements, for many aspects analogous to that of the intersection poset for hyperplane arrangements (see \cite[Section 5]{MoT}). We call the elements of $\mathcal{C}({\widetilde{X}})$ the \emph{layers} of the arrangement. Under our assumptions, the minimal elements of $\mathcal{C}(\widetilde{X})$ are 0-dimensional, hence they are points. We denote by $\mathcal{C}_0(\widetilde{X})$ the set of such layers, which we call the \emph{points} of the arrangement.

For every layer $C$ we define
$${\widetilde{X}}_C\doteq \left\{(\lambda,a)\in {\widetilde{X}} | H_{\lambda,a}\supseteq C\right\}.$$
and
$$X_C\doteq \{\lambda |(\lambda,a)\in {\widetilde{X}}_C\}.$$
The natural surjection ${\widetilde{X}}_C\longrightarrow X_C$ is indeed a bijection, since the condition $(\lambda,a), (\lambda,b)\in X_C$ implies that $\lambda$ is identically equal to $a=b$ on $C$.

\subsection{Primitive vectors}

Given a system of coordinates $(t_1,\ldots,t_n)$ on $T$, for every $$\nu=(\nu_1,\ldots,\nu_n)\in \mathbb{Z}^n$$ we have a map
\begin{align*}
e(\nu):T&\rightarrow\mathbb{C}^*\\
(t_1,\ldots,t_n)&\mapsto{t_1}^{\nu_1}\cdot\ldots\cdot{t_n}^{\nu_n}.
\end{align*}
It is well known that $e$ is an isomorphism between $\mathbb{Z}^n$ and $\Lambda=Hom(T, \mathbb{C}^*)$.

We will assume every $\lambda\in X$ to be \emph{primitive}, i.e. $$\Lambda\cap\langle\lambda\rangle_{\mathbb{C}}=\langle\lambda\rangle_{\mathbb{Z}}.$$
This amounts to require that under the previous isomorphism $\lambda$ is identified with a vector $\nu=(\nu_1,\ldots,\nu_n)\in \mathbb{Z}^n$ such that $GCD(\{\nu_i\})=1$.

\begin{re}\label{pri}
This is not a restrictive assumption; indeed, suppose
$$GCD(\{\nu_i\})=d>1$$ and write $\nu_i'\doteq \nu_i/d$. Then
$${t_1}^{\nu_1}\cdot\ldots\cdot{t_n}^{\nu_n}-a=\left({t_1}^{\nu_1'}\cdot\ldots\cdot{t_n}^{\nu_n'}\right)^d-a=\prod_{i=1}^d \left({t_1}^{\nu_1'}\cdot\ldots\cdot{t_n}^{\nu_n'}-\zeta^i \sqrt[d]{a} \right)$$
where $\zeta$ is a primitive $d-$th root of $1$. Then there is a primitive element $\lambda'$ of $\Lambda$ such that $\lambda=d\lambda'$, and
we can write $H_{\lambda,a}$ as the union of its connected components:
$$H_{\lambda,a}=\bigsqcup_{i=1}^d H_{\lambda',\zeta^i \sqrt[d]{a}}.$$
Then we can replace every pair $(\lambda,a)\in {\widetilde{X}}$ with  all the pairs $(\lambda',\zeta^i a)$. In this way we get a new set $\widetilde{X}'$ which defines the same toric arrangement as $\widetilde{X}$.
\end{re}

\begin{ex}\label{eC2}
Take $T=(\mathbb{C}^*)^2$ with coordinates $(t,s)$ and
$${\widetilde{X}}=\{ (t^2,1),(s^2,1),(ts,1), (ts^{-1},1)\}.$$
Since $t^2-1=(t+1)(t-1)$, the hypersurfaces $H_{t^2}$ and $H_{s^2}$ have two connected components each; $H_{ts}$ and $H_{ts^{-1}}$ are connected, but their intersection is not.

The points of the arrangement are:
$$p_1=(1,1),\: p_2=(-1,-1),\: p_3=(1,-1),\: p_4=(-1,1).$$
Notice that ${\widetilde{X}}_{p_1}={\widetilde{X}}_{p_2}={\widetilde{X}}$, whereas
$${\widetilde{X}}_{p_3}={\widetilde{X}}_{p_4}=\{(t^2,1),(s^2,1)\}.$$

Following Remark \ref{pri}, we can replace ${\widetilde{X}}$ by
$${\widetilde{X}}'=\{(t,1),(t,-1),(s,1),(s,-1),(ts,1), (ts^{-1},1)\}.$$
\end{ex}

\subsection{Construction of the model}\label{mod}

Given a sublattice $\Delta\subset \Lambda$, we define its \emph{completion}
$$\overline{\Delta}\doteq\langle \Delta\rangle_{\mathbb{C}}\cap \Lambda.$$

For every layer $C\in \mathcal{C}({\widetilde{X}})$, we consider the lattice $\lc\doteq \langle X_C\rangle_{\mathbb{Z}}$ and its completion $\lcb$.

\begin{re}\label{lcb}
The elements of $\lcb$ are the characters taking a constant value on $C$.
Indeed, for every $\lambda\in\lcb$, we have that $d \lambda\in\lc$ for some $d>0$. Then by definition $d\lambda$ takes a constant value $a$ on $C$; hence
$$\lambda(t)^{d}=a\;\forall\, t\in C.$$
Since $C$ is connected and the set of $d$th roots of unity is discrete, the continuous map $\lambda$ must be constant.
\end{re}

Now let $\lambda_1,\ldots,\lambda_k$ be an \emph{integral basis} of $\lcb$ (i.e., a basis spanning over $\mathbb{Z}$ the lattice $\lcb$), and let $a_i$ be the constant value assumed by $\lambda_i$ on $C$: then
the ideal $\mathfrak{I}_C$ of the regular functions on $T$ that vanish on $C$ is generated by
$$\left\{\lambda_1- a_1, \ldots, \lambda_k-  a_k\right\}$$
and the normal space to $C$ in $T$ is
$$\mathbf{N}_T(C)\simeq\left(\frac{\mathfrak{I}_C}{\mathfrak{I}_C^2}\right)^*.$$
We denote by $\pc$ its projectified $\pcn$ and by $\varphi_C$ the natural map
\begin{align*}
\varphi_C: T\setminus C & \rightarrow \pc\\
t & \mapsto\: [\lambda_1(t)- a_1,\ldots, \lambda_k(t)- a_k].
\end{align*}
Now let us fix a subset $\g\subseteq \mathcal{C}({\widetilde{X}})$. By collecting the maps $\{\varphi_C,\:C\in \g\}$ and the inclusion $j:\rxt\hookrightarrow T$, we get a map
$$i_\g=j\times\prod_{C\in \g}\varphi_C:\rxt\rightarrow T\times\prod_{C\in \g}\pc$$

We define $\zxg$ as the closure $\overline{i_\g(\rxt)}$ of the image of $\rxt$.

In the next section we will describe the subsets $\g$ that give arise to models with good geometric properties.

\begin{re}

~

\begin{enumerate}\label{pig}
\item If we choose another basis $\lambda'_1, \ldots, \lambda'_k$, we get other generators
$$\{\lambda'_1- a'_1, \ldots, \lambda'_k-  a'_k\}$$
 of the same ideal $\mathfrak{I}_C$, hence another basis of $\mathfrak{I}_C / \mathfrak{I}_C^2$ and then another system of projective coordinates for $\pc$. Then our construction does not depend on such choice.
\item In fact $\zxg$ can be obtained by a sequence of blow-ups along the elements of $\g$, listed in any dimension-increasing order (see \cite{MP}).
\item Since $\prod_{C\in \g}\pc$ is a projective variety, the restriction $\pi: \zxg \rightarrow T$ of the projection on the first factor $T$ is a projective and thus proper map.
  \item Since $i_\g$ is injective, we identify  $\rxt$ with its image  $i_\g(\rxt)$. Such image is closed in $\rxt \times \prod_{C\in \g}\pc$, which is open in $T\times \prod_{C\in \g}\pc$; therefore $\zxg$ contains $\rxt$ as a dense open set, and the restriction of $\pi$ to $\rxt$ is $j$.
\end{enumerate}
\end{re}

\subsection{Hyperplane arrangements and complete sets}

Given a finite set $A\subseteq U$, a \emph{hyperplane arrangement} $\mathcal{H}(A)$ is defined in the dual space
$V=U^*$ by taking the orthogonal hyperplane to each element of
$A$. To every subset $B\subseteq A$ is associated the subspace
$B^\bot$ of $V$ that is the intersection of the corresponding
hyperplanes of $\mathcal{H}(A)$; in other words, $B^\bot$ is the subspace of
vectors that are orthogonal to every element of $B$.
Then we set
$$\mathcal{L}(A)=\{B^\bot, B\subseteq A\}.$$
$\mathcal{L}(A)$ is a geometric lattice, called the \emph{intersection poset} of $\mathcal{H}(A)$. Its elements are called the \emph{flats} of the arrangement.

Given a subset $B\subset A$, we define its completion
$$\overline{B}\doteq\langle B\rangle_{\mathbb{C}}\cap A.$$
We say that $B$ is \emph{complete} in $A$ if $B=\overline{B}$.
In the language of matroid theory, the complete subsets of $A$ are the flats of the associated matroid,
and the flat $\overline{B}$ is the \emph{closure} of the subset $B$.

For every $Q\in\mathcal{L}(A)$, let $\alpha(Q)$ be the set of elements of $A$ which are identically equal to 0 on $Q$; clearly
$$\alpha(Q)^\bot= Q \;\mbox{ and }\;\alpha(B^\bot)=\overline{B}.$$
Hence we have a bijection between $\mathcal{L}(A)$ and the family of complete subsets of $A$.

~

Fix $p\in {\mathcal{C}_0({\widetilde{X}})}$. For every pair $(\lambda,a)\in \widetilde{X_p}$, $\lambda-a\in\mathfrak{I}_p$ defines a vector in ${\mathfrak{I}_p}/{\mathfrak{I}_p^2}$ and hence a hyperplane in its dual, which is the normal space to the point, i.e. the tangent space $T(p)$ to $p$ in $T$. This hyperplane of $T(p)$ is simply the tangent space to the hypersurface $H_{(\lambda,a)}$ in $p$. In this way $X_p$ defines in $T(p)$ a hyperplane arrangement $\mathcal{H}_p$, which is locally isomorphic (in $\underline{0}$) to our toric arrangement (in $p$).
 Then the map
 $$C\mapsto (X_C)^{\bot}$$
 is an inclusion-preserving bijection between layers $C\in \mathcal{C}({\widetilde{X}})$ containing $p$ and flats of $\mathcal{H}_p$.

\begin{re}\label{BiLa}
In particular we see that, for every layer $C$ containing $p$, $X_C=\alpha\big((X_C)^{\bot}\big)$ is a complete subset of $X_p$.
Conversely, for every complete subset $A$ of $X_p$ there is a unique layer $C(A)$ such that $X_{C(A)}=A$ and $p\in C(A)$. Namely, $C(A)$ is the connected component containing $p$ of the subvariety of $T$
$$H_A\doteq\left\{t\in T \; | \; \lambda(t)-\lambda(p)=0\;\forall\,\lambda\in A \right\}.$$
 \end{re}

\section{Combinatorial definitions}

\subsection{Irreducible sets}

Let $B$ be a finite subset of $\Lambda$.
An \emph{integral decomposition} of $B$ is a partition $B=\bigcup_i B_i$ such that
  $$\overline{\langle B\rangle_{\mathbb{Z}}} = \bigoplus_i \overline{\langle B_i\rangle_{\mathbb{Z}}}.$$
A \emph{complex decomposition} of $B$ is a partition $B=\bigcup_i B_i$ such that
  $${\langle B\rangle_{\mathbb{C}}} = \bigoplus_i {\langle B_i\rangle_{\mathbb{C}}}.$$

We say that $B$ is \emph{$\mathbb{Z}-$irreducible} (resp. \emph{$\mathbb{C}-$irreducible}) if it does not have a nontrivial integral (resp. complex) decomposition.

We say that a layer $C\in\mathcal{C}({\widetilde{X}})$ is \emph{$\mathbb{Z}-$irreducible} (resp. \emph{$\mathbb{C}-$irreducible})  if $X_C$ is. We denote by $\mathcal{I}$ (resp. by $\mathcal{I}_\mathbb{C}$) the set of \emph{$\mathbb{Z}-$irreducible} (resp. \emph{$\mathbb{C}-$irreducible}) layers.

\begin{re}
Clearly every integral decomposition is also a complex decomposition, but not conversely: see the example below. Then in general $\mathcal{I}_\mathbb{C}\subsetneq\mathcal{I}$.

In the language of \cite{MP}, $\mathcal{C}({\widetilde{X}})$ is a \emph{conical stratification} on $T$, and $\mathcal{I}_\mathbb{C}$ is the set of the \emph{irreducible strata}. Then a minimal wonderful model can be obtained by blowing up (in any dimension-increasing order) the elements of $\mathcal{I}_\mathbb{C}$. However, in this model the intersections of irreducible components of the normal crossing divisor fail to be connected (see example below). In order to obtain such property (i.e. the last point of Theorem \ref{ncd}), we will blow up all the elements of $\mathcal{I}$.
\end{re}

\begin{ex}
Take $T=(\mathbb{C}^*)^2$ with coordinates $(t,s)$ and
$${\widetilde{X}}=\big\{(ts,1), (ts^{-1},1)\big\}.$$
Then $X$ is identified with the subset $\{(1,1), (1,-1)\}$ of $\mathbb{Z}^2$. Thus $X$
is not $\mathbb{C}-$irreducible, but it is $\mathbb{Z}-$irreducible: indeed
$\mathbb{Z}(1,1)\oplus\mathbb{Z}(1,-1)$ is a sublattice of index $2$ in $\mathbb{Z}^2$.

The hypersurfaces $H_{ts}$ and $H_{ts^{-1}}$ are the irreducible components of a normal crossing divisor; however their intersection consists of two points.
By blowing them up we optain a model whose normal crossing divisor has four irreducible components, pairwise intersecting in a single point (as in the picture below):

\includegraphics[width=120mm]{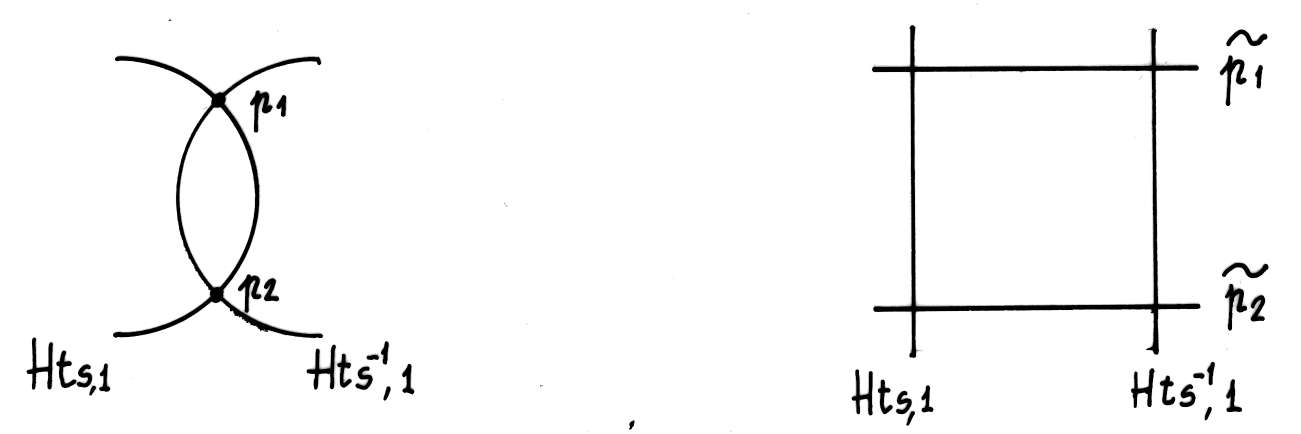}

\end{ex}

We now prove some properties of integral decompositions, which are known (and easier to prove) for complex decompositions (see for instance \cite[Chapter 20.1]{li}).

 From now on we will simply call \emph{decompositions} the integral decompositions, and \emph{irreducible} sets (resp. layers) the $\mathbb{Z}-$irreducible sets (resp. layers).

\begin{lem}
Let $B=B_1\cup B_2$ be a decomposition and $D\subset B$ be an irreducible subset. Then $D\subseteq B_1$ or $D\subseteq B_2$.
\end{lem}

\begin{proof}
Set $D_1\doteq D\cap B_1$ and $D_2\doteq D\cap B_2$. We must prove that $D=D_1\cup D_2$ is a decomposition; then the irreducibility of $D$ implies that $D_1$ or $D_2$ is empty.
We first notice that
$${\langle D\rangle_{\mathbb{Z}}}={\langle D_1\rangle_{\mathbb{Z}}}\oplus{\langle D_2\rangle_{\mathbb{Z}}}$$
since $${\langle D_1\rangle_{\mathbb{Z}}}\cap{\langle D_2\rangle_{\mathbb{Z}}}\subseteq{\langle B_1\rangle_{\mathbb{Z}}}\cap{\langle B_2\rangle_{\mathbb{Z}}}\subseteq \overline{\langle B_1\rangle_{\mathbb{Z}}}\cap\overline{\langle B_2\rangle_{\mathbb{Z}}}=\{\underline{0}\}.$$

Then take any $\lambda\in \overline{\langle D\rangle_{\mathbb{Z}}}$. For some positive integer $m$ we have that $m\lambda\in{\langle D\rangle_{\mathbb{Z}}}$ and then can be written uniquely as $m\lambda=\mu_1+\mu_2$, with $\mu_1\in {\langle D_1\rangle_{\mathbb{Z}}}$ and $\mu_2\in {\langle D_2\rangle_{\mathbb{Z}}}$. Moreover, since
$$\lambda\in \overline{\langle B\rangle_{\mathbb{Z}}}=\overline{\langle B_1\rangle_{\mathbb{Z}}}\oplus\overline{\langle B_2\rangle_{\mathbb{Z}}}$$
 $\lambda$ can be expressed uniquely as $\lambda=\gamma_1+\gamma_2$, with $\gamma_1\in \overline{\langle B_1\rangle_{\mathbb{Z}}}$ and $\gamma_2\in \overline{\langle B_2\rangle_{\mathbb{Z}}}$. Then $m\lambda=m\gamma_1+m\gamma_2=\mu_1+\mu_2$ implies $\mu_1=m\gamma_1$ and $\mu_2=m\gamma_2$, hence $\gamma_1\in \overline{\langle D_1\rangle_{\mathbb{Z}}}$ and $\gamma_2\in \overline{\langle D_2\rangle_{\mathbb{Z}}}$. Thus
$$\overline{\langle D\rangle_{\mathbb{Z}}}=\overline{\langle D_1\rangle_{\mathbb{Z}}}\oplus\overline{\langle D_2\rangle_{\mathbb{Z}}}.$$
\end{proof}

\begin{lem}
Every subset $B$ has a decomposition $B=\bigcup B_i$ into irreducible subsets $B_i$. This decomposition is unique up to the order.
\end{lem}

\begin{proof}
The existence is clear by induction. Now let $B=\bigcup B'_j$ be another decomposition into irreducible subsets. By the previous lemma every $B_i$ is contained in some $B'_j$ and vice versa. Then these factors are the same up to the order.
\end{proof}

\subsection{Building sets and nested sets of layers}

We now recall some general definitions given in \cite{DPw} and \cite[Chapter 20.1]{li}, adapting them to our situation.

Let $A$ be a finite subset of $\Lambda$.
A family $\g^*$ of subsets of $A$ is a \emph{building set} if every complete subset $B$ of $A$ is decomposed by the maximal elements $B_i$ of $\g^*$ contained in $B$. Then we say that $B=\bigcup_i B_i$ is the decomposition of $B$ in $\g^*$ or that the $B_i$s are the $\g^*-$factors of $B$.

A subset $\s^*$ of $\g^*$ is a \emph{$\g^*-$nested set} if given any $B_1,\ldots, B_r\in \s^*$ mutually incomparable,
      $$B\doteq B_1\cup\ldots\cup B_r$$
       is a complete set in $A$ with its decomposition in $\g^*$.

By \cite{MP}, an equivalent definition is the following. A \emph{flag} $\mathcal{F}^*$ is a sequence $A_1\subset\dots\subset A_k$ of subsets of $A$.
A set $\s^*=\{B_1,\dots, B_s\}$ is  $\g^*-$\emph{nested} if there is a flag $\mathcal{F}^*$ such that all the elements of $\s^*$ are $\g^*-$factors of elements of $\mathcal{F}^*$.

The family $\mathcal{I}^*$ of all irreducible subsets of $A$ is clearly a building set.
In particular, we call nested sets the $\mathcal{I}^*-$nested sets. Then a \emph{nested set} is a family $\s^*$ of irreducible subsets such that for every $B_1,\ldots, B_r\in \s^*$ mutually incomparable,
      $$B\doteq B_1\cup\ldots\cup B_r$$
       is a complete set in $A$ with its decomposition into irreducible subsets.

~



%
%

Now let $p\in \mathcal{C}_0(\widetilde{X})$  be a point of the arrangement, and let $C$ be any layer containing $p$.
Let $\g^*$ be a building set in $X_p$, and let $X_C=\bigcup_{i}X_i$ be the decomposition of $X_C$ in $\g^*$.
By Remark \ref{BiLa}, there is a unique layer $C_i\doteq C(X_i)$ containing $C$ and such that $X_{C_i}=X_i$.
We call the $C_i$s the \emph{$\g-$factors} of $C$; clearly $C=\cap C_i$.

Then we can associate to every building set $\g^*$ in $A$ a \emph{building set of layers} $\g$ defined as the set of all the $\g-$factors of all the elements of $\mathcal{C}({\widetilde{X}})$. In particular for $\g^*=\mathcal{I}^*$ we get that the set $\mathcal{I}$ of all irreducible layers is a building set.

A \emph{flag} $\mathcal{F}$ of layers is a sequence $D_1\subset\dots\subset D_k.$
A set of layers
$$\s=\{C_1,\dots, C_s\}$$
is \emph{$\g-$nested} if there is a flag $\mathcal{F}$ such that all the elements of $\s$ are $\g-$factors of elements of $\mathcal{F}$.
We say that $\s$ is a \emph{nested set of layers} if it is $\mathcal{I}-$nested, i.e. if there is a flag $\mathcal{F}$ such that all the elements of $\s$ are irreducible factors of elements of $\mathcal{F}$.

We call the minimal element of the flag (with respect to inclusion) the \emph{center} of $\s$. This is a well defined layer by the following Lemma:

\begin{lem}
Let $\s$ be a $\g-$nested set. Then
$$C(\s)\doteq\bigcap_{C\in\s}C$$
is connected (and then is a layer).
\end{lem}

\begin{proof}
Let $M(\s)$ be the set of minimal elements of $\s$, with respect to inclusion. Clearly
$$C(\s)=\bigcap_{C\in M(\s)}C.$$
The elements of $M(\s)$ are pairwise incomparable, hence
$$\overline{\Lambda_{C(\s)}}=\sum_{C\in S}\overline{\Lambda_C}=\bigoplus_{C\in M(\s)}\overline{\Lambda_C}.$$
Let us choose an integral basis $\underline{b}_C$ for each of the lattices $\overline{\Lambda_C}$, $C\in M(\mathcal S)$. Then
$$\underline{b}=\bigcup_{C\in M(\s)}\underline{b}_C$$
 is an integral basis for $\overline{\Lambda_{C(\mathcal S)}}$. For any $\lambda\in \overline{\Lambda_C}$, $\lambda$ takes a constant value $a_{\lambda}$ on $C$ by Remark \ref{lcb}. It follows that the elements $\lambda-a_{\lambda}$, $\lambda\in\underline{b}$ generate the ideal of definition of $C(\mathcal S)$, which is clearly irreducible since $\underline{b}$ is a basis of a split direct summand in $\Lambda$.
\end{proof}
\begin{re}
Notice that our proof clearly implies that the intersection $$C(\mathcal S)=\bigcap_{C\in M(\mathcal S)} C$$ is transversal.
\end{re}

A $\g-$nested set of layers is \emph{maximal}  if it is not contained in a larger one; this happens if and only if $\s$ contains all the irreducible $\g-$factors of a maximal flag.
In this case the center of $\s$ is a point $p=p(\s)$. We denote by $\mathfrak{M}$ the set of all maximal $\g-$nested set of layers of $\mathcal{C}({\widetilde{X}})$ and by $\mathfrak{M}_p$ the set of those having center $p$. Then we have the partition
$$\mathfrak{M}=\bigsqcup_{p\in {\mathcal{C}_0({\widetilde{X}})}}\mathfrak{M}_p.$$

The following fact is clear from the definitions (and from Remark \ref{BiLa}):

\begin{lem}\label{bss}
If $\s=\{C_1,\ldots, C_s\}\in\mathfrak{M}_p$  is a maximal $\g-$nested set of layers of center $p$, then $$\s^*\doteq\{X_{C_1},\ldots, X_{C_s}\}$$ is a maximal $\g^*-$nested set in $X_p$.

Conversely, given a maximal $\g^*-$nested set $\widehat{\s}$ in $X_p$, there is a unique $\s\in\mathfrak{M}_p$ such that $\s^*=\widehat{\s}$; namely
$$\s\doteq \left\{C(A_i),\; A_i\in \widehat{\s} \right\}.$$
\end{lem}

In particular $|\s|=|\s^*|=n$, the rank of $X$ (see \cite[Theor 20.9]{li}).

Finally we prove an elementary result that we will use frequently in the next sections. Take $\s\in\mathfrak{M}_p$.
\begin{lem}\label{suc}

~

\begin{enumerate}
  \item Let $C\in \mathcal{I}$ and $p\in C$. Then there is an element ${\overline{C}}\in\s$ which is the maximum among all the elements of $\s$ contained in $C$; we call it the \emph{$\s-$core} of $C$.
  \item Let $C$ be an element of $\s$ which is not minimal in it. Then there is an element $s(C)\in\s$ which is the maximum among all the elements of $\s$ properly contained in $C$; we call it the \emph{successor} of $C$.
\end{enumerate}
\end{lem}
\begin{proof}
The proof is the same for both statements. Let $C'$ and $C''$ be two elements of $\s$ which are contained (or, for the second statement, properly contained) in C.
Then $X_C\subset X_{C'}\cap X_{C''}$; hence $X_{C'}\cup X_{C''}$ is not a decomposition. Since $X_{C'}$ and $X_{C''}$ are in the $\g^*-$nested set $\s^*$, they must be comparable; then also $C'$ and $C''$ are.
\end{proof}

\subsection{Adapted bases}

Given a $\g-$nested set $\s$, we say that an integral basis $\underline{b}\doteq\{\lambda_1\,\ldots,\lambda_n\}$ for the lattice $\Lambda$ is \emph{adapted to $\s$} if for every $C\in \s$, $\underline{b}\cap \overline{\Lambda_C}$ is an integral basis for $\overline{\Lambda_C}$.

\begin{lem}\label{ada}
 There exists an integral basis $\underline{b}^\s$ for $\Lambda$ adapted to $\s$.
 \end{lem}
\begin{proof}
Let us define
$$\Lambda_\s\doteq {\sum_{D\in \s}\overline {  \Lambda_D }}.$$
Notice that
$$\Lambda_{\mathcal S}=\bigoplus_{C\in M(\mathcal S)}\overline {  \Lambda_D }$$
where $M(\mathcal S)$ is the set of minimal (and hence pairwise incomparable) elements of $\mathcal S$. Then by definition $\Lambda_{\mathcal S}=\overline {\Lambda_{\mathcal S}}$.
We will prove, by induction on the cardinality of $\s$, that there is a basis of $\Lambda_\s$ adapted to $\s$.
Then our claim follows: indeed, since the lattice $\Lambda_\s$ either coincide with $\Lambda$ or is a split direct summand of it, the basis of $\Lambda_\s$ can be completed to a basis of $\Lambda$.

If $\s$ contains only one element $C$, the statement is trivial since $\Lambda_\s=\overline{\Lambda_C}$ and every basis of this lattice is adapted to $\s$.

Otherwise, take a minimal $C\in \s$, and set $\s'=\s\setminus \{C\}.$
Since $\s'$ is $\g-$nested, by inductive hypothesis the lattice
$$\Lambda_{\s'}= {\sum_{D\in \s'}\overline {  \Lambda_D }}$$
has an integral basis adapted to $\mathcal S'$.
Since $\Lambda_{\mathcal S'}=\overline {\Lambda_{\mathcal S'}}$ we can complete the chosen basis of $\Lambda_{\mathcal S'}$ to an integral basis $\underline b$ of $\Lambda_{\mathcal S}$ using elements of $\overline {  \Lambda_C }$. We claim that this basis is adapted to $\mathcal S$. Let us take $D$ in $\mathcal S$. If $D\neq C$ there is nothing to prove. Then assume $D=C$. In this case we know that
$$\Lambda_{\mathcal S}=\overline {  \Lambda_C }\oplus \bigoplus_{D\in M(\mathcal S)\setminus \{C\}} \overline {  \Lambda_D }.$$
By construction, every element in $\underline b$ either lies in $\overline {  \Lambda_C }$ or in
$$\bigoplus_{D\in M(\mathcal S)\setminus \{C\}} \overline {\Lambda_D }.$$
Then every $\lambda\in\lcb$ is in the span of $\underline{b}\cap\lcb$,
 proving our claim.
\end{proof}

To every maximal $\g-$nested set of layers $\s\in\mathfrak{M}_p$ we associate a function
$$p_\s:\Lambda\longrightarrow \s$$
in the following way.
For every $\lambda\in \Lambda$ we set $a\doteq \lambda(p)$, and we define
$p_{\s}(\lambda)$ as the maximum element of $\s$ on which $\lambda$ is identically equal to $a$. This is well defined by Lemma \ref{suc}: indeed $p_{\s}(\lambda)= \overline{H_{(\lambda, a)}}$. This function has the following properties:

\begin{lem}\label{psl}

~

\begin{enumerate}
  \item For every $C\in\g$ there exists $\lambda\in X_C$ such that $p_\s(\lambda)=\overline{C}$.
  \item The restriction of $p_\s$ to an adapted basis $\underline{b}$ is a bijection.
\end{enumerate}
\end{lem}

\begin{proof}
For every $C\in\g$, let $M(C)$ be the (possibly empty) set of the elements of $\s$ properly containing $C$ and minimal with this property. Such elements are pairwise incomparable, hence $$\bigcup_{D\in M(C)}X_D$$ is a decomposition. Since $X_C\supset X_D$ for every $D\in M(C)$,
$$X_C\supset \bigcup_{D\in M(C)}X_D$$
and this inclusion is proper, because $X_C\in\g^*$. Then there exists
$$\lambda\in X_C\setminus \bigcup_{D\in M(C)}X_D.$$
By definition $p_{\s}(\lambda)=\overline{C}$, then the first statement is proved.

Now assume $C\in\s$, and let $\underline{b}$ be an adapted basis to $\s$: then by definition $\underline{b}\cap \overline{\Lambda_C}$ is a basis for $\overline{\Lambda_C}$ and
$$\bigsqcup_{D\in M(C)}\left(\underline{b}\cap \overline{\Lambda_D}\right)\mbox{ is a basis for }\bigoplus_{D\in M(C)} \overline{\Lambda_D}.$$
Since $C\in\g$, we have that $$\overline{\Lambda_C}\supsetneq \bigoplus_{D\in M(C)} \overline{\Lambda_D}.$$
Then there exists
$$\lambda\in\left(\underline{b}\cap \overline{\Lambda_C}\right)\setminus \bigsqcup_{D\in M(C)}\left(\underline{b}\cap \overline{\Lambda_D}\right).$$
Clearly $p_{\s}(\lambda)=C$. Then we proved that the restriction of $p_\s$ to $\underline{b}$ is surjective; therefore it is bijective, since $|\underline{b}|=n=|\s|$.
\end{proof}

\bigskip

\begin{re}
From now on we will assume for simplicity $\g=\ir$, and then we will focus on the model $\zxt\doteq\mathbf{Z}_{{\widetilde{X}},\ir}$ defined as the closure of the image of the map
$$i_\ir=j\times\prod_{C\in \ir}\varphi_C:\rxt\rightarrow T\times\prod_{C\in \ir}\pc.$$
However, all the results in this paper may be extended to the case of an arbitrary building set $\g$.
\end{re}

\section{Open sets and smoothness}

\subsection{Definition of the open sets}

To every $\s\in\mathfrak{M}_p$ we associate a nonlinear change of coordinates $f_\mathcal{S}$ and an open set $\vsp$ defined as follows.

Let us take a basis of $\Lambda$ adapted to $\s$, and denote it by $$\underline{b}^\s=\left(\lambda_C\right)_{C\in\s}$$ where $\lambda_C\doteq p_{\s}^{-1}(C)$. Set $a_C\doteq\lambda_C(p)$.
Since $\underline{b}^\s$ is integral, $\left(\lambda_C-a_C\right)_{C\in\s}$ is a system of coordinates on $T$.

Consider $\mathbb{C}^n$ with coordinates $\underline{z}^\s=(z_C)_{C\in \s}$, and its open set $$\widetilde{U_\s}\doteq \left\{(z_C)\in \mathbb{C}^n | \prod_{D\subseteq C} z_D\neq - a_C \:\forall C\in\s\right\}.$$

Define a map $f_\s :\widetilde{U_\s}\rightarrow T$ in the given coordinates as
$$\lambda_C\left(f_\s(\underline{z}^\s)\right)=\left(\prod_{D\subseteq C}z_D\right)+ a_C$$
or equivalently as the nonlinear change of coordinates
\begin{equation}\label{nlc}
    \lambda_C-a_C=\prod_{D\subseteq C}z_D.
\end{equation}
Then  $f_\s(\underline{0})=p$.

Notice that on the open set of $T$ where $\lambda_C-a_C\neq 0\forall C\in \s$,
the map $f_\s$ can be inverted by the following formula:
\begin{equation}\label{inv}
z_C=\begin{cases}
\lambda_{{C}}- a_{{C}} & \mbox{ , if ${C}$ is minimal in $\s$}\\
\frac{\lambda_{{C}}\,- a_{{C}}}{\lambda_{s({C})}- a_{s({C})}} & \mbox{ ,  otherwise}
\end{cases}\end{equation}
where $s({C})$ is the successor defined in Lemma \ref{suc}.

Let us define the open set of $T$
$$T_p\doteq T\setminus \bigcup_{p\notin C} C$$
and set ${U_\s} \doteq {f_\s}^{-1}(T_p)$. We denote again by $f_\s$ the restriction ${U_\s}\rightarrow T_p$.

Now take any $\lambda\in \Lambda$; set $a\doteq \lambda(p)$ and $C\doteq p_\s(\lambda)$.

Since $\underline{b}^\s$ is adapted to $\s$, an integral basis for  $\lcb$ is given by
$$\underline{b}^\s\cap \overline{\Lambda_{{C}}}=\left\{\lambda_D, D\supseteq {C}   \right\}.$$
In particular 
$\lambda$ can be expressed in this basis, and since $p_\s(\lambda)=C$, $\lambda$
does not lie in the span of $\left\{\lambda_D, D\supsetneq C   \right\}$: then
$$\lambda=m_C\lambda_C+\sum_{D\supsetneq C}m_D\lambda_D$$
for some integers $m_D$ and a nonzero integer $m_C$.
The previous identity, considered as an equality of regular functions on $T$, can be written as
$$\lambda=\lambda_C^{m_C}\prod_{D\supsetneq C}\lambda_D^{m_D}.$$
Then we have:
\begin{equation}\label{agt}
\lambda-a=\left(\lambda_C^{m_C}\prod_{D\supsetneq C}\lambda_D^{m_D}-a_C^{m_C}\prod_{D\supsetneq C}\lambda_D^{m_D}\right)+\left(a_C^{m_C}\prod_{D\supsetneq C}\lambda_D^{m_D}-a\right)
\end{equation}
and we can write the first summand as
$$\prod_{D\supsetneq C}\lambda_D^{m_D}\left(\lambda_C^{m_C}-a_C^{m_C}\right)=\beta_C(\lambda_C-\alpha_C)$$
where
$$\beta_C\doteq\prod_{D\supsetneq C}\lambda_D^{m_D}\prod_{\zeta^{m_C}=1,\zeta\neq 1}\left(\lambda_C-\zeta a_C\right)$$
is a regular function on $T$ which is invertible on $C$. Working in the same way on the second summand of Formula (\ref{agt}) we see that,
for some regular functions $\{\beta_D, D\in \s\}$,
$$\lambda- a=\beta_{ {C}}\,(\lambda_{ {C}}- a_{ {C}})+\sum_{D\supsetneq  {C}}\beta_D(\lambda_D- a_D).$$
By operating the change of coordinates (\ref{nlc}), we get:
\begin{equation}\label{pl1}
\lambda- a=
\left(\beta_{ {C}} \prod_{E\subseteq  {C}}z_E+\sum_{D\supsetneq  {C}}\beta_D\prod_{E\subseteq D}z_E\right)
=\left(\prod_{E\subseteq  {C}}z_E\right)\cdot p_\lambda(\underline{z}^\s)
\end{equation}
where we set
\begin{equation*}
    p_\lambda(\underline{z}^\s)\doteq \beta_{ {C}}+\sum_{D\supsetneq  {C}}{\beta_D}\prod_{D\supseteq E\supsetneq  {C}}z_E.
\end{equation*}

We define $\vsp$ 
as the open set 
of ${U_\s}$ where
$$\prod_{\lambda\in X_p}p_\lambda(\underline{z}^\s)\neq 0.$$

Let us remark that $\underline{0}\in\vsp$, since for every $\lambda\in X_p$ we have that $p_\lambda(\underline{0})=\beta_C(p) \neq 0$.
 Furthermore in $\vsp$, for every $\lambda\in X_p$, we have the equality of regular functions
\begin{equation}\label{div}
\prod_{E\subseteq p_\s(\lambda)}z_E = \frac{\lambda- a}{p_\lambda(\underline{z}^\s)}.
\end{equation}

\subsection{Properties of the open sets}

Let us define the open set of $\vsp$
$$\vsp^0\doteq \{\underline{z}\in\vsp \:|\: z_C\neq 0\forall C\in \s \}.$$

We denote by $A_\s$ the open set of $T$ given by $f_\s(\vsp)\cap\rxt$. We remark that by Formula (\ref{div}) $f_\s^{-1}(A_\s)=\vsp^0$ and the restriction of $f_\s$ to $\vsp^0$ maps it into $A_\s$. By composing this map with the inclusion $A_\s\hookrightarrow\rxt$ and with the application $\phi_C:\rxt\rightarrow \pc$ defined in Section \ref{mod}, we get a map
$$\psi_C:\vsp^0\longrightarrow \pc.$$

\begin{lem}
For every $C\in \ir$ and $\s\in \mathfrak{M}_p$, the map $\psi_C$
extends uniquely to a map
$$\widetilde{\psi_C}:\vsp\rightarrow \pc.$$
\end{lem}

\begin{proof}
Let $p$ be the center of $\s$.
If $C$ does not contain $p$ the statement is clear: indeed since $\vsp\subset U_\s$, for every $u\in\vsp$ we have that $t\doteq f_\s(u)\notin C$ so that
for at least one index $j$, $\lambda_j(t)\neq a_j$.
Then the projective coordinate $\lambda_j(t)- a_j$ of $\pc$ is nonzero.

Then assume $p\in C$, and let $\overline{C}$ be its $\s-$core (see Lemma \ref{suc}).
By the first part of Lemma \ref{psl}, there exists $\lambda_1\in X_C$ such that $p_\s(\lambda_1)=\overline{C}$. Since we assumed (Remark \ref{pri}) every element of $X_C$ to be primitive, we can complete $\{\lambda_1\}$ to an integral basis $\{\lambda_1,\ldots,\lambda_k\}$ of $\overline{\Lambda_C}$. Then if we set $a_i\doteq \lambda_i(p)$, we have that $$[\lambda_1- a_1,\ldots, \lambda_k- a_k]$$ is a system of projective coordinates for $\pc$.

Since $\underline{b}^\s$ is adapted to $\s$, an integral basis for $\overline{\Lambda_{\overline{C}}}$ is given by
$$\underline{b}^\s\cap \overline{\Lambda_{\overline{C}}}=\left\{\lambda_D, D\supseteq \overline{C}   \right\}.$$
In particular every $\lambda_i\in \lcb\subseteq \overline{\Lambda_{\overline{C}}}$ can be expressed in this basis, and since $p_\s(\lambda_1)=\overline{C}$, $\lambda_1$ does not lie in the span of $\left\{\lambda_D, D\supsetneq \overline{C}   \right\}.$

After making the nonlinear change of coordinates (\ref{nlc}) as in Formula (\ref{pl1}), we can divide every projective coordinate by $\prod_{E\subseteq \overline{C}}z_E$; in this way we get that the map $\psi_C:\vsp^0\longrightarrow \pc$ is given by
$$\underline{z}\mapsto\left[p_{\lambda_1}(\underline{z}),\; p_{\lambda_2}(\underline{z})\prod_{\overline{C}\subsetneq E\subseteq D_2}z_E,\;\ldots,\; p_{\lambda_k}(\underline{z})\prod_{\overline{C}\subsetneq E\subseteq D_k}z_E \right]$$
where we set $D_i\doteq p_\s(\lambda_i)$. Since by definition $p_{\lambda_1}(\underline{z})\neq 0$ for $\underline{z}\in\vsp$, this map extends to $\vsp$. Moreover its image is contained in an affine open set of $\pc$.

Finally the uniqueness of the extension is clear since by its very definition $\vsp^0$ is dense in $\vsp$.
\end{proof}

By applying the lemma above to all the layers $C\in \ir$, we get that for every $\s\in \mathfrak{M}_p$ the inclusion  
$\vsp^0\hookrightarrow \zxt$
extends uniquely to a map
$$j_\s:\vsp\rightarrow \zxt.$$

\begin{lem}
The map $j_\s$ is an embedding into a smooth open set.
\end{lem}

\begin{proof}
In order to prove that $j_\s$ is an embedding, it suffices to see that every coordinate $z_C$ on $\vsp$ can be written as the composition of $j_\s$ and a function on $j_\s(\vsp)$.
Then take $C\in \s$. If $C$ is not minimal, let $D=s(C)$ be the successor of $C$. Since $\underline{b}^\s$ is adapted to $\s$, on $\mathbb{P}_D$ we have the projective coordinates
$$\left[\lambda_E- a_E\right]_{E\in\s, E \supseteq D}$$
and by the proof of the previous lemma $\vsp$ maps into the affine subset where $\lambda_D- a_D\neq 0$.
Then we can read the coordinate $z_C$ in $\mathbb{P}_D$ by Formula (\ref{inv}):
$$z_C=\frac{\lambda_C- a_C}{\lambda_D- a_D}.$$
If on the other hand $C$ is minimal in $\s$, then $z_C=\lambda_C- a_C$.

In this way all the coordinates $z_C$ can be recovered by the projection of $j_\s(\vsp)\subset\zxt$ on $T$ or on some $\mathbb{P}_D$; hence our map is an embedding.
Moreover, since $(z_C)_{C\in\s}$ is a system of coordinates on $j_\s(\vsp)$, in every point the differential of $j_\s$ has rank $|\s|=n$.
Then $j_\s(\vsp)$ is smooth.
\end{proof}

\begin{re}
By abuse of notation, from now on we will write $\vsp$ for $j_\s(\vsp)$, identifying this set with its isomorphic image in $\zxt$.
\end{re}

\subsection{Smoothness of the model}
Let us define $$\yxt\doteq \bigcup_{\s\in\mathfrak{M}}\vsp .$$
In this section we prove that $\yxt=\zxt$, and hence $\zxt$ is smooth. The main step is the following lemma, which tells that every curve in $\rxt$ that "has limit" in $T$, "has limit" in $\yxt$.
Let $D_\varepsilon\doteq \left\{s\in \mathbb{C} \; | \;  |s|<\varepsilon\right\}$.
\begin{lem}
Let $f: D_\varepsilon \rightarrow T$ be a curve such that $f(D_\varepsilon\setminus \{0\})\subseteq \rxt$.

Then $f$ lifts to a curve in $\yxt$.
\end{lem}

\begin{proof}
Given such a $f$, let $C_f\in \mathcal{C}({\widetilde{X}})$ be the smallest layer containing $f(0)$, and let $p\in \mathcal{C}_0({\widetilde{X}})$ be a point contained in $C_f$.
For every $\lambda\in X_p$, we have that locally, near $s=0$, we can write
$$\lambda(f(s))- a=s^{n_\lambda}q_\lambda(s)$$
with $ a=\lambda(p)$, $n_\lambda\geq 0$ and $q_\lambda(0)\neq 0$.

For every integer $h\geq 0$, let us define
$$A_h\doteq \{ \lambda\in X_p | n_\lambda\geq h\}.$$
Notice that $A_0=X_p$ and $A_{h+1}\subseteq A_h$; by taking all the irreducible factors of the elements of this flag we get a nested set in $X_p$. Let us complete it to a maximal nested set $\s^*$; by Lemma \ref{bss}, to $\s^*$ is naturally associated a maximal nested set of layers $\s\in \mathfrak{M}_p$.

We claim that for a such $\s$, the curve $f: D_\varepsilon\setminus \{0\}\rightarrow \rxt$ extends to a map $f:D_\varepsilon\rightarrow \vsp$.

First notice that $f(0)\in T_p$: indeed for every layer $D$ containing $f(0)$ we have that $C_f\subseteq D$ by minimality and then $p\in D$.
Then we have to prove that:
\begin{enumerate}
       \item $z_C\big(f(s)\big)$ is defined in $0$ for every $C\in \s$;
       \item $p_\lambda\big(f(0)\big)\neq 0$ for every $\lambda\in X_p$.
\end{enumerate}

Take $C\in \s$; if $C$ is minimal in $\s$ then $z_C(f(s))=\lambda_C(f(s))- a_C$ and there is nothing to prove. Otherwise,
let $D=s(C)$ be the successor of $C$. Then by Formula (\ref{inv})
$$z_C(f(s))=\frac{\lambda_C(f(s))- a_C}{\lambda_D(f(s))- a_D}=s^{n_{\lambda_C}-n_{\lambda_D}}\frac{q_{\lambda_C}(s)}{q_{\lambda_D}(s)}$$
and $n_{\lambda_C}\geq n_{\lambda_D}$ by the definition of $\s$, so $z_C$ is well defined in $0$.

As for the second claim, given any $\lambda\in X_p$ set $C\doteq p_\s(\lambda)$ and take the vector $\lambda_C$ of the adapted basis $\underline{b}^\s$.

Then by definition of $\s$, $n_{\lambda}= n_{\lambda_{C}}$, and by Formulae (\ref{nlc}) and (\ref{pl1}) we have
$$p_\lambda=\frac{\lambda- a}{\lambda_C- a_C}.$$

Therefore
$$p_\lambda(f(0))=\frac{\lambda\big(f(0)\big)- a}{\lambda_{{C}}\big(f(0)\big)- a_{{C}}}=
\frac{q_{\lambda}(0)}{q_{\lambda_{{C}}}(0)}\neq 0.$$
\end{proof}

\begin{te}\label{smo}
$\yxt=\zxt$. In particular $\zxt$ is smooth.
\end{te}
\begin{proof}
By the well known \emph{valuative criterion for properness} (see for instance \cite{Har}),
the previous lemma amounts to say that the map
$$\pi |_{\yxt}:\yxt\rightarrow T$$ is proper.
Since also the projection $$T\times\prod_{C\in \mathcal{I}}\pc\rightarrow T$$ is proper, the embedding
$$\yxt\rightarrow T\times\prod_{C\in \mathcal{I}}\pc$$
is proper as well; therefore its image is closed, and thus it coincides with $\zxt$.

Therefore $\zxt$ is smooth, since it is the union of smooth open sets.
\end{proof}

\section{The normal crossing divisor}

\subsection{Technical lemmas}

For every $C\in \mathcal{I}$, let us define a divisor $\mathbf{D}_C\subset \zxt$ as follows. Take a $\s\in \mathfrak{M}$ such that $C\in\s$. In the open set $\vsp$ take the divisor of equation $z_C=0$; let $\mathbf{D}_C$ be the closure of this divisor in $\zxt$. The following lemma implies that $\mathbf{D}_C$ does not depend on the choice of $\s$, and yields the theorem below, which describes the geometry of $\zxt\setminus\rxt$.

\begin{lem}\label{coo}
Take any two maximal nested sets of layers $\s\in \mathfrak{M}_p$ and $\mathcal{Q}\in \mathfrak{M}_q$. Let $\{\zs_C, C\in \s\}$ and $\{\zt_C, C\in \mathcal{Q}\}$ be the corresponding sets of coordinates on $\vsp$ and $\vtq$.

Then for every $C\in\s$:
\begin{enumerate}
  \item if $C\in\s\setminus\mathcal{Q}$, $\zs_C$ is invertible as a function on $\vsp\cap\vtq$;
  \item if $C\in \s\cap \mathcal{Q}$, $\zs_C / \zt_C$ is regular and invertible as a function on $\vsp\cap\vtq$.
\end{enumerate}
\end{lem}
\begin{proof}
If $q\notin C$, then $C\in\s\setminus\mathcal{Q}$, and the (first) statement is proved as follows. Take $x\in\zxt$ such that $z^\s_C(x)=0$: then by Formula (\ref{nlc}) $\pi(x)\in C$, where $\pi: \zxt\rightarrow T$ is the projection defined in Remark \ref{pig}. Therefore $\pi(x)\notin T_q$, hence $x\notin\vtq$, proving the claim.

Therefore we can assume $q\in C$ and proceed by induction as in the proof of \cite[Lemma 20.39]{li}:
\begin{itemize}
  \item First let us assume $C$ to be a minimal element in $\mathcal{I}$; then necessarily $C\in \s\cap \mathcal{Q}$. We recall that $\zs_C=\lambda^{\s}_C- a^{\s}_C$; set
  $$D\doteq p_{\mathcal{Q}}(\lambda^{\s}_C)\supseteq C.$$
Then for some function $a$
$$\zs_C=a\prod_{E\in \mathcal{Q},D\supseteq E}\zt_E=a\zt_C \prod_{E\in \mathcal{Q}, D\supseteq E, E\neq C}\zt_E.$$
In the same way $\zt_C=\lambda^{\mathcal{Q}}_C- a^{\mathcal{Q}}_C$, and if we set
  $$D'\doteq p_{\s}(\lambda^{\mathcal{Q}}_C)\supseteq C$$
we get
$$\zt_C=a'\prod_{F\in \s,D'\supseteq F}\zt_F=a'\zs_C \prod_{F\in \s, D'\supseteq F, F\neq C}\zs_F.$$
for some function $a'$.
Since both $D$ and $D'$ contain $C$, by substituting we get:
$$\zs_C=\zs_C \; a \,a' \prod_{E\in \mathcal{Q}, D\supseteq E, E\neq C}\zt_E\prod_{F\in \s, D'\supseteq F, F\neq C}\zs_F.$$
Therefore
$$a a' \prod_{E\in \mathcal{Q},D\supseteq E, E\neq C}\zt_E\prod_{F\in \s,D'\supseteq F, F\neq C}\zs_F=1$$
and hence
$$\frac{\zs_C}{\zt_C}=a \prod_{E\in \mathcal{Q},D\supseteq E, E\neq C}\zt_E$$
is invertible, as claimed.

  \item Now let us take any $C\in\s$. By induction, we can assume that our claims are true for every $D\subsetneq C$, $D\in\s\cup \mathcal{Q}$
  (if $D\in\mathcal{Q}\setminus\s$, by symmetry $\zt_D$ is assumed to be invertible on $\vsp\cap\vtq$).

 Let $D={\overline{C}}\in\mathcal{Q}$ be the $\mathcal{Q}-$core of $C$. Take $\lambda\in X_C$ such that $p_\mathcal{Q}(\lambda)=D$, and set $G\doteq p_\s (\lambda)$. Then
 $G\supseteq C$ and $\lambda$ takes on $D$ and on $G$ the same constant value $ a\doteq\lambda(p)$.
 Notice that $D$ is the $\mathcal{Q}-$core of $G$.

 Then for some invertible $b, b'$
 $$\lambda- a= b \prod_{E\in \mathcal{Q}, D\supseteq E}\zt_E = b'\prod_{F\in \s, G\supseteq F}\zs_F.$$
 Hence
 \begin{equation}\label{prz}
    1=b^{^{-1}}\, b' \prod_{F\in \s\setminus \mathcal{Q}, G\supseteq F}\zs_F \prod_{E\in \mathcal{Q}\setminus\s, D\supseteq E}{\zt_{E}}^{-1} \prod_{F\in \s\cap \mathcal{Q}, D\supseteq F}\zs_F{\zt_{F}}^{-1}.
 \end{equation}

We can now prove the first claim. If $C\notin \mathcal{Q}$ then $D\subsetneq C$. Then all the factors in equation (\ref{prz}) are regular: those of type
$$\zs_F, F\in \s\setminus \mathcal{Q}, G\supseteq F$$
obviously, the others by inductive assumption, since they involve elements properly contained in $C$. Since $\zs_C$ appears as one of the factors in (\ref{prz}) it is invertible.

In the same way if $C\in \mathcal{Q}$, and then $D=C$, all the factors in (\ref{prz}) but (eventually) $\zs_C\, {\zt_C}^{-1}$ are regular; then also $\zs_C\, {\zt_C}^{-1}$ must be regular and invertible.
\end{itemize}
\end{proof}

\begin{lem}\label{dsg}
Let be $C\in \mathcal{I}$.
\begin{enumerate}
  \item The divisor $\mathbf{D}_C$ is well defined.
  \item If $C\notin\s$, then $\mathbf{D}_C\cap \vsp=\emptyset$.
\end{enumerate}
\end{lem}
\begin{proof}
\begin{enumerate}
  \item Let $\s, \mathcal{Q}$ be two maximal nested set of layers containing $C$. Then by the second point of Lemma \ref{coo}, $\zs_C$ and $\zt_C$ have the same zeros in $\vsp\cap\vtq$, which is an open dense set in $\vsp$ and in $\vtq$. Then the closures of the two divisors coincide.
  \item Let $\mathcal{Q}$ be a maximal nested set of layers containing $C$. Then by the first point of Lemma \ref{coo}, $\zt_C$ is invertible as a function on $\vsp\cap\vtq$. Therefore the divisor of $\vtq$ defined by $\zt_C=0$ is contained in $\zxt\setminus\vsp$. Since this set is closed, it also contains $\mathbf{D}_C$ which is the closure of the divisor.
\end{enumerate}
\end{proof}

\subsection{The main theorem}

Now let us define
$$\mathbf{D}=\bigcup_{C\in \mathcal{I}}\mathbf{D}_C.$$
The geometry of the divisor $\mathbf{D}$ is described by the following theorem.

\begin{te}\label{ncd}

~

\begin{enumerate}
  \item $\zxt\setminus \mathbf{D}=\rxt$.
  \item $\mathbf{D}$ is a normal crossing divisor whose irreducible components are the divisors $\mathbf{D}_C, C\in \mathcal{I}$.
  \item Let be $\mathcal{N}\subseteq \mathcal{I}$, and
  $$\dn\doteq \bigcap_{C\in \mathcal{N}}\mathbf{D}_C.$$
    Then $\dn\neq\emptyset$ if and only if $\mathcal{N}$ is nested.
  \item If $\mathcal{N}$ is nested, $\dn$ is smooth and irreducible.
\end{enumerate}
\end{te}
\begin{proof}
By Theorem \ref{smo}, we can check each statement on every open set $\vsp, \s\in \mathfrak{M}$.

Then the first claim, by the second part of Lemma \ref{dsg}, amounts to note that
$$\left(\zxt\setminus \mathbf{D}\right)\cap\vsp=\vsp\setminus\bigcup_{C\in\s}\left(\mathbf{D}_C\cap\vsp\right)=\vsp^0=\rxt\cap\vsp.$$
(for the definition of $\vsp^0$ see the beginning of Section 4.2).

The second statement is obvious since
$$\mathbf{D}\cap\vsp=\bigcup_{C\in\s}\left(\mathbf{D}_C\cap\vsp\right)=\left\{\underline{z}\in\vsp | z_C= 0 \mbox{ for some } C\in\s\right\}$$
is by definition a normal crossing divisor in $\vsp$.

For the third statement, note that if $\mathcal{N}$ is not nested it is not contained in any maximal nested set of layers; then for every $\s\in \mathfrak{M}$, $\dn\cap\vsp=\emptyset$ by the second part of Lemma \ref{dsg}. On the other hand, if $\mathcal{N}$ is nested it can be completed to some $\s\in \mathfrak{M}$, and
$$\dn\cap\vsp=\left\{\underline{z}\in\vsp | z_C= 0 \forall C\in\mathcal{N}\right\}$$
which is clearly nonempty, smooth and irreducible. Since
$$\dn=\bigcup_{\s\supseteq\mathcal{N}}\left(\dn\cap\vsp\right)$$
also the last statement follows.

\end{proof}


\end{document}